\documentclass[12pt]{amsart}
\usepackage{amsmath,amssymb}
\newtheorem{theorem}{Theorem}
\newtheorem{proposition}[theorem]{Proposition}
\def\B{\Bbb B}
\def\C{\Bbb C}
\def\CC{\mathcal C}
\def\K{\mathcal K}
\def\h{h}
\def\ve{\varepsilon}
\def\ds{\displaystyle}

\title[The squeezing function near \h-extendible boundary points]
{Behavior of the squeezing function near \h-extendible boundary points}

\author{Nikolai Nikolov}

\address{Institute of Mathematics and Informatics\\Bulgarian Academy
of Sciences\\ Acad. G. Bonchev 8, 1113 Sofia, Bulgaria\newline
\indent Faculty of Information Sciences\\
State University of Library Studies and Information Technologies\\
Shipchenski prohod 69A, 1574 Sofia,
Bulgaria}\email{nik@math.bas.bg}



\begin{document}

\begin{abstract} It is shown that if the squeezing function tends to one
at an h-extendible boundary point of a $\CC^\infty$-smooth, bounded
pseudoconvex domain, then the point is strictly pseudoconvex.
\end{abstract}

\maketitle

Denote by $\B_n$ the unit ball in $\C^n.$ Let $M$ be an $n$-dimensional complex manifold,
and $z\in M.$ For any holomorphic embedding $f:M\to\B_n$ with $f(z)=0,$ set
$$s_M(f,z)=\sup\{r>0:r\B_n\subset f(M)\}.$$
The squeezing function of $M$ is defined by $\ds s_M(z)=\sup_f s_M(f,z)$
if such $f$'s exist, and $\ds s_M(z)=0$ otherwise.

Many properties and applications of the squeezing function  have been explored by various authors,
see e.g. \cite{DGZ,FW,JK} and the references therein.

It was shown in \cite{DGZ} that if $D$ is a $\CC^2$-smooth strictly pseudoconvex domain
in $\C^n,$ then
\begin{equation}\label{1}
\lim_{z\to\partial D}s_D(z)=1.
\end{equation}

A. Zimmer \cite{Z1} proved the converse if $D$ is a $\CC^\infty$-smooth, bounded convex
domain; namely, if \eqref{1} holds, then $D$ is necessarily strictly pseudoconvex. Recently,
he extended this result to the $\CC^{2,\alpha}$-smooth case \cite{Z2}.

On the other hand, J.E. Forn{\ae}ss and E.F. Wold \cite{FW} provided an example showing that
$\CC^2$-smoothness is not enough. They also asked if Zimmer's result holds for $\CC^\infty$-smooth,
bounded pseudoconvex domains.

S. Joo and K.-T. Kim \cite{JK} gave an affirmative answer for domains of finite type in $\C^2.$

This can be extended to larger class of domains by using different arguments.
Recall that a $\CC^\infty$-smooth boundary point $a$ of finite type of a domain $D$ in $\C^n$
is said to be h-extendible \cite{Yu1,Yu2} (or semiregular \cite{DH})
if $D$ is pseudoconvex near $a,$ and Catlin's and D'Angelo multitypes of $a$ coincide.

For example, $a$ is extendible if the Levi form at $a$ has a corank at most one \cite{Yu2}, or
$D$ is linearly convexifiable near $a$ \cite{Con}. In particular, h-extendibility takes
place in the strictly pseudoconvex, two-dimensional finite type, and convex finite type cases.

\begin{theorem}\label{main} Let $a$ be an h-extendible boundary point of a $\CC^\infty$-smooth,
bounded pseudoconvex domain $D$ in $\C^n.$ If $s_D(a_j)\to 1$ for a nontangential sequence $a_j\to a,$
then $a$ is a strictly pseudoconvex point.
\end{theorem}

Nontangentiality means that $\ds\liminf_{j\to\infty}\frac{d_D(a_j)}{|a_j-a|}>0,$ where $d_D$ is the distance
to $\partial D.$
\smallskip

Before proving Theorem \ref{main}, we need some preparation.

Denote by $\mu=(m_1,m_2,\dots,m_n)$ Catlin's multitype of $a$
($m_1=1$ and $m_2\le\dots\le m_n$ are even numbers). By \cite{DH,Yu1,Yu2}, there exists
a local change of variables $w=\Phi(z)$ near $a$ such that $\Phi(a)=0,$ $J\Phi(a)=1,$
$$r(\Phi^{-1}(w))=\mbox{Re}(w_1)+P(w')+o(\sigma(w)),$$
where $r$ is the signed distance to $\partial D,$ $\ds\sigma(z)=\sum_{j=1}^n|w_j|^{m_j}$ and
$P$ is a $1/\mu$-homogeneous polynomial without pluriharmonic terms. Moreover, the so-called model domain
$$E=\{w\in\C^n:\mbox{Re}(w_1)+P(w')<0\}$$
(which depends on $\Phi$) is of finite type.

In \cite{Yu2,Nik}, the nontangential boundary behavior of the Kobayashi-Royden and Carath\'eodory-Reiffen
metrics of $D$ near $a$ are expressed in terms of $r,$ $\Phi,$ and the respective metrics of $E_{\Phi}$
at its interior point $e=(-1,0').$ Obvious modifications in the proofs of these results allows to obtain
similar results for the Kobayashi-Eisenman and Carath\'eodory-Eisenman volumes of $D:$
$$\K_D(u)=\inf\{|Jf(0)|^{-1}:f\in\mathcal O(\B^n,D),f(0)=u\},$$
$$\CC_D(u)=\sup\{|Jf(u)|:f\in\mathcal O(D,\B^n), f(u)=0\}.$$

\begin{proposition}\label{vol} Let $a$ be an h-extendible boundary point of a domain $D$ in $\C^n.$
Let $\mu$ be Catlin's multitype of $a$ and $\ds m=\sum_{j=1}^n\frac1{m_j}.$ Then
\begin{equation}\label{k}
\K_D(a_j)(d_D(a_j))^m\to \K_E(e)
\end{equation}
for any nontangential sequence $a_j\to a.$

If, in addition, $D$ is $\CC^\infty$-smooth, bounded pseudoconvex, then
\begin{equation}
\label{c}\CC_D(a_j)(d_D(a_j))^m\to\CC_E(e).
\end{equation}
\end{proposition}

Since $E$ is hyperbolic with respect to the Carath\'eo\-dory-Reiffen
metric \cite{Nik}, it is easy to see $\CC_E>0.$ So, the limits above are positive.
\smallskip

\noindent{\it Sketch of the proof of Proposition \ref{vol}.} Let $\ve>0$ and  
$$E_{\pm\ve}=\{w\in\C^n:\mbox{Re}(w_1)+P(w')\pm\ve\sigma(w)<0\}.$$
There exists a neighborhood $U_\ve$ of $a$ such that
$$E_{+\ve}\cap V_\ve\subset\Phi(D\cap U_\ve)\subset E_{-\ve}\cap V_\ve,$$ where $V_\ve=\Phi(U_\ve).$ 
Since $a\in\partial D$ is a local holomorphic peak point \cite{DH,Yu1},
the localization $\ds\frac{\K_{D\cap U_\ve}(a_j)}{K_D(a_j)}\to 1$ holds.
On the other hand, $E_{\pm\ve}$ are taut domains if $\ve\le\ve_0$ \cite{Yu2}; in particular,
$\K_{E_{\pm\ve}}$ are continuous functions. Let $-c_j+id_j$ be the first coordinate of $b_j$
($c_j,d_j\in\Bbb R$). Since $d_j/c_j$ is a bounded sequence, 
it suffices to show \eqref{k}, when $d_j/c_j\to s.$ Set $b_j=\Phi(a_j)$ and 
$\pi_j(w)=(w_1c_j^{-1/m_1},\dots,w_n c_j^{-1/m_n}).$ Note that $\pi_j(b_j)\to e_s:=(-1+is,0').$ Now, applying 
the scaling of coordinates $\pi_j$ and using normal family arguments, we obtain that 
$\K_{E_{\pm\ve}\cap V_\ve}(b_j)c_j^m\to\K_{E_{\pm\ve}}(e_s).$ Finally, following \cite[Theorem 2.1]{Yu2}, 
one can prove that $\K_{E_{\pm\ve}}(e_s)\to \K_E(e)$ as $\ve\to 0.$ These facts, together with 
$\ds\frac{c_j}{d_D(a_j)}\to 1,$ imply \eqref{k}.

The proof of \eqref{c} follows by similar but more delicate arguments
as in \cite{Nik}.\qed
\smallskip

\noindent{\it Proof of Theorem \ref{main}.} By \cite{DGZ}, one has that
$$(s_D(a_j))^n\K_D(a_j)\le\CC_D(a_j)\le\K_D(a_j).$$
It follows by Proposition \ref{vol} that
$$\CC_E(e)=\K_E(e).$$
Since $\B_n$ and $E$ are taut domains \cite{Yu1,Yu2},
there exist extremal functions for $\CC_E(e)$
and $\K_E(e).$ Then the Carath\'eodory-Cartan-Kaup-Wu theorem
implies that $E$ and $\B_n$ are biholomorphic. Since
$E$ is a model domain of finite type, the main result in \cite{CP}
shows that $m_2=\dots=m_n=2,$ that is, $a$ is a strictly pseudoconvex point.\qed

\end{document}